\input amstex
\documentstyle{amsppt}
\NoRunningHeads
\pagewidth{14cm}
\pageheight{19.1cm}
\magnification=1200
\loadmsbm
\loadbold

\define\ol{\overline}
\define\vvp{\varphi}

\NoBlackBoxes

\dedicatory Dedicated to Bruce Allison on the occasion of his 60th birthday
\enddedicatory

\author Yoji Yoshii
\endauthor

\affil Department of Mathematics\\ North Dakota State University \\  Fargo, ND, 58105-5075  USA \\
yoji.yoshii\@ndsu.nodak.edu
\endaffil

\topmatter
\title Cayley Polynomials
\endtitle

\abstract 
We consider a polynomial version of the Cayley numbers.
Namely,
we define the ring of Cayley polynomials
in terms of generators and relations in
the category of alternative algebras. The ring turns out to be an octonion algebra 
over an ordinary polynomial ring.  Also, a localization 
(a ring of quotients) of 
the ring of Cayley polynomials gives another description
of an octonion torus.
Finally, we find a subalgebra of a prime nondegenerate alternative algebra
so that the subalgebra is an octonion algebra over the center.

\endabstract

\footnotetext""{2000 Mathematics Subject Classification 17D05.}
\endtopmatter

\document

\head
Introduction
\endhead

Nonassociative analogues of Laurent polynomials naturally appeared
in the classification of extended affine Lie algebras and Lie tori.
These Lie algebras are a natural generalization of affine Kac-Moody Lie algebras
(see [AABGP], [N], [Y2]).
As the affine Kac-Moody Lie algebras are coordinatized by the ring
of Laurent polynomials in one variable,
extended affine Lie algebras or Lie tori are coordinatized by
nonassociative analogues of Laurent polynomials in several variables.
Those Lie algebras have types classified by finite irreducible root systems,
and the coordinate algebras depend on the types.
In particular, such Lie algebras of type $\text A_2$
are coordinatized by alternative algebras,
and an alternative analogue of Laurent polynomials
(which is not associative)
was found in [BGKN].  
The coordinate algebra is called an {\it octonion torus}.
It turns out that
the coordinate algebras of
extended affine Lie algebras or Lie tori of type $\text A_2$, $\text C_3$ or $\text F_4$,
which are not associative,
are exactly octonion tori
(see also [AG], [BY], [Y1], [Y3]).

An octonion torus 
(an octonion $n$-torus)
is defined by a Cayley-Dickson process over a ring of Laurent polynomials.
More precisely, it is obtained by
the Cayley-Dickson process three times over $F[z_1^{\pm 1},\ldots,z_n^{\pm 1}]$ with $n\geq 3$,
where $F$ is a field of characteristic $\neq 2$,
taking the structure constants $z_1$, $z_2$ and $z_3$, i.e.,
in the standard notation for the Cayley-Dickson process (see \S 1),
$$(F[z_1^{\pm 1},\ldots,z_n^{\pm 1}],z_1,z_2,z_3).$$
To study the algebra structure, 
it is enough to consider the case $n=3$, the octonion $3$-torus
$(F[z_1^{\pm 1},z_2^{\pm 1},z_3^{\pm 1}],z_1,z_2,z_3)$,
which is also called the {\it Cayley torus}.
Our goal is to find a simple presentation of the Cayley torus, 
or essentially, a simple presentation of its subalgebra
$$D=(F[z_1,z_2,z_3],z_1,z_2,z_3).$$
%generated by the monomials of nonnegative powers of the basic generators in the Cayley torus.
The algebra $D$ also appears as a certain subalgebra of a free alternative algebra
generated by more than three elements,
which was discovered by Dorofeev (see Remark 3.2).

We define the algebra $F_C[t_1,t_2,t_3]$ over $F$ 
by the following relations
$$t_2t_1=-t_1t_2,
\quad
t_3t_1=-t_1t_3,
\quad t_3t_2=-t_2t_3
\quad\text{and}\quad
(t_1t_2)t_3=-t_1(t_2t_3)
$$
in the category of alternative algebras,
and call it the ring of {\it Cayley polynomials} or a {\it universal octonion algebra}
(since it covers all the octonion algebras over $F$).
Then
we show that $F_C[t_1,t_2,t_3]$ is isomorphic to $D$.
In particular,
$F_C[t_1,t_2,t_3]$ is an octonion algebra over the center
$F[t_1^2,t_2^2,t_3^2]$
(ordinary commutative associative polynomials in three variables $t_1^2$, $t_2^2$ and $t_3^2$).
Also,
the Cayley torus can be viewed as the ring of quotients of $F_C[t_1,t_2,t_3]$
by the monomials of the center $F[t_1^2,t_2^2,t_3^2]$.
As corollaries, we obtain a simple presentation of the Cayley torus
and also a presentation of any octonion algebra over $F$.
Moreover,
the base field $F$ can be generalized to a ring $\Phi$ of scalars
containing $1/2$,
and so we will set up the notions above
over $\Phi$.

Finally, we will discuss about {\it Cayley-Dickson rings}.
Such a ring embeds into an octonion algebra over a field.
However, the ring itself is not necessarily an octonion algebra
in general.
We note that our ring $F_C[t_1,t_2,t_3]$ 
or an octonion torus is a Cayley-Dickson ring
and also an octonion algebra.
We will see in Proposition 4.2 that there exists
a subring $B$ of a Cayley-Dickson ring $R$ 
%containing $1/2$
(or of a prime nondegenerate alternative algebra $R$ over $\Phi$)
so that $B$ is an octonion algebra
over the center of $R$,
and the central closure $\ol B$ coincides with the central closure $\ol R$.

\medskip

We thank Professor Bruce Allison and Erhard Neher for several suggestions.

\medskip

Throughout the paper let $\Phi$ be a unital commutative
associative ring containing $1/2$.
Also, all algebras are assumed to be unital.

\head
\S1 Cayley-Dickson Process
\endhead

We recall the Cayley-Dickson Process over a ring $\Phi$ of scalars
(see [M] for detail).
For an algebra $B$ over $\Phi$,
we assume that $B$ is {\it faithful},
i.e., for all $\alpha\in\Phi$,
$\alpha 1=0\Longrightarrow \alpha =0$.
Let $*$ be a {\it scalar involution} of $B$ over $\Phi$, i.e.,
an anti-automorphism of period 2
with $bb^{*}\in\Phi 1$.
Let $\mu\in \Phi$
be a {\it cancellable scalar},
i.e., 
$\mu b=0$ for some $b\in B$ $\Longrightarrow$ $b=0$.
The Cayley-Dickson algebra (or process) over $\Phi$
with structure constant $\mu$
constructed from $B=(B,*)$ is a new algebra $B\oplus B$
with product 
$(a,b)(c,d)=(ac+\mu db^{*},a^{*}d+cb)$
for $a,b,c,d\in B$. 
Letting $v=(0,1)$ we can write this algebra as $B+ vB$
with multiplication
$$(a+vb)(c+vd)=(ac+\mu db^{*})+v(a^{*}d+cb).
\tag1.1
$$
We call $v$ the {\it basic generator}.
Note that $v^2=\mu$.
The algebra $B+ vB$ also has the new involution $*$, defined by
$(a+vb)^{*}=a^*-vb$, which is scalar.
So one can continue the process,
and we write, for example,
$(B,\mu,\nu)$ instead of  $((B,\mu),\nu))$.
Note that 
$B$ is faithful $\Rightarrow$ 
$(B, \mu)$ is faithful, and
$\nu$ is cancellable for $B$ $\Leftrightarrow$ 
$\nu$ is cancellable for $(B, \mu)$.

Let $\mu_{1}, \mu_{2},\mu_{3}$ be any cancellable scalars of $\Phi$.
The Cayley-Dickson process twice starting from $\Phi$
with trivial involution, say
$(\Phi, \mu_{1}, \mu_{2})$, is called
a {\it quaternion algebra},
which is a 4-dimensional free $\Phi$-module
and an associative  but not commutative algebra,
and three times, say
$(\Phi, \mu_{1}, \mu_{2},\mu_{3})$, 
is called an {\it octonion algebra},
which is an 8-dimensional free $\Phi$-module
and an alternative but not associative algebra.
Note that quaternion algebras and octonion algebras are central,
and if $\Phi$ is a field, they are simple.

\proclaim{Lemma 1.2}
Let $v_1$, $v_2$ and $v_3$ be the basic generators in each step of 
an octonion algebra $(\Phi, \mu_{1}, \mu_{2},\mu_{3})$
so that $v_1^2=\mu_1$, $v_2^2=\mu_2$ and $v_3^2=\mu_3$.
Then $v_2v_1=-v_1v_2$,
$v_3v_1=-v_1v_3$,
$v_3v_2=-v_2v_3$
and
$(v_1v_2)v_3=-v_1(v_2v_3)$.
\endproclaim
\demo{Proof}
One can easily check these identities from (1.1).
\qed
\enddemo

We will consider quaternion algebras and octonion algebras over
various rings of scalars, not necessarily $\Phi$
in the following sections.

\head
\S2 Hamilton Polynomials
\endhead

The associative algebra over $\Phi$
with generators $t_1$ and $t_2$
and the relation $t_1t_2=-t_2t_1$
is called the ring of {\it Hamilton polynomials}
or a {\it universal quaternion algebra},
denoted $\Phi _H[t_1,t_2]$.
Note that the center of $\Phi _H[t_1,t_2]$ is equal to $\Phi [t_1^2,t_2^2]$
(the ordinary commutative associative polynomials 
over $\Phi$ in two variables $t_1^2$ and $t_2^2$),
and $\Phi _H[t_1,t_2]$ is a quaternion algebra over $\Phi [t_1^2,t_2^2]$,
i.e., $(\Phi [t_1^2,t_2^2],t_1^2,t_2^2)$
using the notation in \S 1.
(Consider the base ring as $\Phi [t_1^2,t_2^2]$ instead of $\Phi$.
Then $t_1^2$ and $t_2^2$ are cancellabe elements of $\Phi [t_1^2,t_2^2]$.)
Also,
it is clear that any quaternion algebra over $\Phi $ is a homomorphic image of $\Phi _H[t_1,t_2]$.

Note that 
the multiplicative subset $S:=\{t_1^rt_2^s\}_{r,s\in 2\Bbb N}$
of the center $\Phi [t_1^2,t_2^2]$
does not contain zero divisors of $\Phi _H[t_1,t_2]$,
and so
one can construct the ring of quotients $S^{-1}\Phi _H[t_1,t_2]$
 (see e.g. [SSSZ, p.185]).
Then $S^{-1}\Phi _H[t_1,t_2]$ is still a quaternion algebra, i.e.,
$$S^{-1}\Phi _H[t_1,t_2]=(\Phi [t_1^{\pm 2},t_2^{\pm 2}], t_1^2, t_2^2).$$
Give the degrees $(1,0)$, $(0,1)$, $(-1,0)$ and $(0,-1)$
for $t_1$, $t_2$, $t_1^{-1}$ and $t_2^{-1}$,
respectively.
(Note that $t_1$ and $t_2$
are invertible in $S^{-1}\Phi _H[t_1,t_2]$
and $t_1^{-1}=t_1^{-2}t_1$ and $t_2^{-1}=t_2^{-2}t_2$.)
Then $S^{-1}\Phi _H[t_1,t_2]$ becomes a $\Bbb Z^2$-graded algebra, called the 
{\it quaternion $2$-torus} or the {\it Hamilton torus}.
Note that $\Phi _H[t_1,t_2]$ embeds into $S^{-1}\Phi _H[t_1,t_2]$.

If $A$ is an associative algebra over $\Phi $
generated by invertible elements $a$ and $b$,
and they satisfy $ab=-ba$,
then $A$ is a homomorphic image of
$S^{-1}\Phi _H[t_1,t_2]$ via $t_1\mapsto a$ and $t_2\mapsto b$,
using the universal property of the ring of quotients.
Also, the associative algebra $L$ over $\Phi $ with generators $t_1^{\pm 1}$ and $t_2^{\pm 1}$
and relations $t_1t_1^{-1}=t_2t_2^{-1}=1$ has a natural $\Bbb Z^2$-grading
as above.
So there is a natural graded homomorphism from $S^{-1}\Phi _H[t_1,t_2]$
onto $L/(t_1t_2+t_2t_1)$.
On the other hand,
since $S^{-1}\Phi _H[t_1,t_2]$ has the relations defining $L/(t_1t_2+t_2t_1)$,
there is a natural graded homomorphism from $L/(t_1t_2+t_2t_1)$
onto $S^{-1}\Phi _H[t_1,t_2]$.
Hence they are graded isomorphisms.
Thus 
the Hamilton torus $S^{-1}\Phi _H[t_1,t_2]$
has a presentation in the category of associative algebra;
generators $t_1^{\pm 1}$ and $t_2^{\pm 1}$
with relations 
$t_1^{-1}t_1=t_2t_2^{-1}=1$ and $t_1t_2=-t_2t_1$.
Because of the presentation,
it is reasonable to write
$$S^{-1}\Phi _H[t_1,t_2]=\Phi _H[t_1^{\pm 1},t_2^{\pm 1}].$$
Also, a {\it quaternion $n$-torus} ($n\geq 2$) is defined as
$$
\Phi _H[t_1^{\pm 1},\ldots, t_n^{\pm 1}]:=
\Phi _H[t_1^{\pm 1},t_2^{\pm 1}]\otimes_\Phi \Phi[t_3^{\pm 1},\ldots,t_n^{\pm 1}],
$$
where $\Phi[t_3^{\pm 1},\ldots,t_n^{\pm 1}]$ is 
the ordinary Laurent polynomial algebra over $\Phi$
in $(n-2)$-variables. 

\medskip

The following proposition is well-known in ring theory.
One can prove it in the same way as in Theorem 3.6 for
octonion algebras.

\proclaim{Proposition 2.1}
Any quaternion algebra over $\Phi $, say 
$(\Phi ,\mu_1,\mu_2)$ for cancellable scalars $\mu_1,\mu_2$ of $\Phi $,
 is isomorphic to $\Phi _H[t_1,t_2]/(t_1^2-\mu_1, t_2^2-\mu_2)$.
Hence $(\Phi ,\mu_1,\mu_2)$ has a presentation
in the category of associative algebras;
generators $t_1$ and $t_2$
with relations 
$t_1^2=\mu_1$, $t_2^2=\mu_2$ and $t_1t_2=-t_2t_1$.

In particular, if $\Phi $ is a field, $A$ is an associative algebra over $\Phi $
generated by $a_1$ and $a_2$,
and they satisfy $a_1a_2=-a_2a_1$,
$a_1^2=\mu_1$ and $a_2^2=\mu_2$,
then $A$ is isomorphic to
$(\Phi ,\mu_1,\mu_2)$.

\endproclaim

\head
\S3 Cayley Polynomials
\endhead

We will use
the commutator $[a,b]=ab-ba$ and the associator $(a,b,c)=(ab)c-a(bc)$
in the subsequent claims.
Alternative algebras are defined by two idendtiites:
$(a,a,b)=0=(b,a,a)$.
We have the alternative law
$(a,b,c)=-(b,a,c)=(b,c,a)$, etc.,
and the flexible law
$(a,b,a)=0$, and so
we can omit the parentheses for $(ab)a=a(ba)$.
We will use the middle Moufang identity
$(ab)(ca)=a(bc)a$ in Proposition 3.1.
Recall that the center of an alternative algebra $A$ is defined
as $\{z\in A\mid [z,a]=(z,a,b)=0\ \text{for all $a,b\in A$}\}$.

\enskip

The  
alternative algebra  over $\Phi $ with generators $t_1$, $t_2$, $t_3$
and the {\it Cayley relations} 
$$t_2t_1=-t_1t_2,
\quad
t_3t_1=-t_1t_3,
\quad t_3t_2=-t_2t_3
\quad\text{and}\quad
(t_1t_2)t_3=-t_1(t_2t_3)
\tag C
$$
is called the ring of {\it Cayley polynomials}
or a {\it universal octonion algebra},
denoted $\Phi _C[t_1,t_2,t_3]$.
Note that any octonion algebra over $\Phi $ is a homomorphic image of $\Phi _C[t_1,t_2,t_3]$
by Lemma 1.2.

Let
$Z$
be the center of $\Phi _C[t_1,t_2,t_3]$.
Our main goal is to show that $\Phi _C[t_1,t_2,t_3]$ is an octonion algebra over $Z$.
\proclaim{Claim 1}
$\Phi _C[t_1,t_2,t_3]$ has the identities
$(t_it_j)t_k=-t_i(t_jt_k)$
(anti-associativity) and 
$(t_it_j)t_k=-t_k(t_it_j)$
(anti-commutativity)
for any distinct $i,j,k\in\{1,2,3\}$.
\endproclaim
\demo{Proof}
By the anti-commutativity in (C), it suffices to show
three identities for the anti-commutativity $(t_it_j)t_k=-t_k(t_it_j)$,
say $k=1,2,3$.
However, we need to prove five identities for the anti-associativity
$(t_it_j)t_k=-t_i(t_jt_k)$.

By the alternative law, we have
$(t_1,t_2,t_3)=-(t_2,t_1,t_3)$.
So
$$
\align
(t_2t_1)t_3+t_2(t_1t_3)
&=(t_2t_1)t_3-(t_2,t_1,t_3)+(t_2t_1)t_3\\
&=(t_1,t_2,t_3)-2(t_1t_2)t_3\quad\text{(since $t_2t_1=-t_1t_2$)}\\
&=0\quad\text{(since $(t_1t_2)t_3=-t_1(t_2t_3)$)}.
\endalign
$$
Hence,
$$(t_2t_1)t_3=-t_2(t_1t_3)
\quad\text{(anti-associativity)}.
\tag a1$$
Also, from
$(t_1,t_2,t_3)=-(t_1,t_3,t_2)$,
we have
$(t_1t_2)t_3-t_1(t_2t_3)=-(t_1t_3)t_2+t_1(t_3t_2)$.
Since $t_3t_2=-t_2t_3$,
we get 
$$(t_1t_2)t_3=-(t_1t_3)t_2.
\tag1$$
For the rest of argument, we will use the identities in (C)
without mentioning.
By (1), we have
$(t_1t_3)t_2+t_1(t_3t_2)
=-(t_1t_2)t_3+t_1(t_3t_2)=0$.
Hence,
$$
(t_1t_3)t_2=-t_1(t_3t_2)
\quad\text{(anti-associativity)}.
\tag a2$$
By (a1) and (1), we get
$$t_2(t_1t_3)=-(t_1t_3)t_2
\quad\text{(anti-commutativity for $k=2$)}.
\tag2$$
Also, 
$(t_3t_1)t_2+t_3(t_1t_2)
=-(t_1,t_2,t_3)+2(t_3t_1)t_2=-(t_1,t_2,t_3)-2(t_1t_3)t_2
=-(t_1,t_2,t_3)+2(t_1t_2)t_3=0$
by (1).
Hence,
$$(t_3t_1)t_2=-t_3(t_1t_2)
\quad\text{(anti-associativity)}.
\tag a3$$
By (1) and (a3), we get
$$t_3(t_1t_2)=-(t_1t_2)t_3
\quad\text{(anti-commutativity for $k=3$)}.
\tag3$$
Now, 
$(t_3t_2)t_1+t_3(t_2t_1)
=(t_3,t_2,t_1)+2t_3(t_2t_1)
=-(t_1,t_2,t_3)-2t_3(t_1t_2)
=-(t_1,t_2,t_3)+2(t_1t_2)t_3=0$
by (3).
Hence,
$$(t_3t_2)t_1=-t_3(t_2t_1)
\quad\text{(anti-associativity)}.
\tag a4$$
By (a4), the main involution 
on the free alternative algebra
$\Phi \langle t_1,t_2,t_3\rangle$,
i.e., the involution determined by
$t_1\mapsto t_1$, $t_2\mapsto t_2$ and $t_3\mapsto t_3$,
preserves the relations (C).
So applying for the induced involution to (a2),
we get 
$$
(t_2t_3)t_1=-t_2(t_3t_1)
\quad\text{(anti-associativity)}.
\tag a5$$
Finally, 
$$
\align
(t_3t_2)t_1+t_1(t_3t_2)
&=-t_3(t_2t_1)+t_1(t_3t_2)
\quad\text{by (a4)}\\
&=t_3(t_1t_2)-t_1(t_2t_3)
\\
&=-(t_1t_2)t_3-t_1(t_2t_3)
\quad\text{by (3)}\\
&=0.
\endalign
$$
Hence,
$$(t_3t_2)t_1=-t_1(t_3t_2)
\quad\text{(anti-commutativity for $k=1$)}.
\qed
$$
\enddemo

\proclaim{Claim 2}
$t_1^2,t_2^2,t_3^2\in Z$.
\endproclaim
\demo{Proof}
We have $[t_1^2, t_i]=(t_1^2,t_i,t_j)=0$ for all $i,j\in\{1,2,3\}$.
Indeed, except for the cases $(i,j)=(2,3)$ and $(3,2)$ in the second identity,
this follows from (C) and Artin's Theorem, that is, any subalgebra generated by 2 elements is associative.
For the cases $(i,j)=(2,3)$ and $(3,2)$,  
one can use the identity
$(a^2, b, c)=(a, ab+ba,c)$
[SSSZ, (17), p.36]
for any alternative algebra.
Hence, $(t_1^2,t_2,t_3)=(t_1^2,t_3,t_2)=0$ by (C).
Thus, $t_1^2$ is central for a generating set of $\Phi _C[t_1,t_2,t_3]$,
and by the theorem of Bruck and Kleinfeld [SSSZ, Lemma 16, p.289], 
we obtain $t_1^2\in Z$.
By the symmetry of our relations (C) with Claim 1, we also obtain $t_2^2, t_3^2\in Z$.
\qed
\enddemo
\proclaim{Claim 3}
$(t_1t_2)^2=-t_1^2t_2^2$,  $(t_1t_3)^2=-t_1^2t_3^2$, $(t_2t_3)^2=-t_2^2t_3^2$
and $(t_1(t_2t_3))^2=t_1^2t_2^2t_3^2$.
In particular, they are all in the center $Z$.
\endproclaim
\demo{Proof}
Again, by Artin's Theorem, we have $(uv)^2=uvuv$ (no parentheses are needed).
So, for example,
$(t_1t_2)^2=-t_1^2t_2^2$, or for the last one,
$(t_1(t_2t_3))^2=t_1(t_2t_3)t_1(t_2t_3)=-(t_2t_3)t_1^2(t_2t_3)$ (by Claim 1)
$=-t_1^2(t_2t_3)^2=t_1^2t_2^2t_3^2$.
\qed
\enddemo

We can give  a natural $\Bbb N^3$-grading to $\Phi _C[t_1,t_2,t_3]$,
defining $\deg t_1=(1,0,0)$, $\deg t_2=(0,1,0)$
and $\deg t_3=(0,0,1)$.
This is possible because
$\Phi _C[t_1,t_2,t_3]$ is defined by the homogeneous relations (C).

\proclaim{Proposition 3.1}
Let $t=t_{i_1}\cdots t_{i_k}$ be an element 
of degree $(\ell,m,n)$ in $\Phi _C[t_1,t_2,t_3]$,
omitting various parentheses, where $t_{i_j}=t_1$, $t_2$ or $t_3$,
so that the total degree of $t$ is $\ell+m+n=k$.
Then $t=\pm (t_1^\ell t_2^m)t_3^n$.
\endproclaim
\demo{Proof}
It is clear for $k=1$, 2 or 3 by Claim 1.
Suppose that the total degree $k>3$.
Then there exists at least one of the following parts in $t$:
(i) $a(bc)$, (ii) $(ab)c$ or (iii) $(ab)(cd)$
for some $a,b,c,d\in\{t_1,t_2,t_3\}$.
For (i), if two of $a,b,c$ are the same, then $a(bc)=\pm t_p^2t_q$
for some $p,q\in\{1,2,3\}$
by Claim 1 and 2.
Hence, $t=\pm t't_p^2$ for some $t'$ 
which is a product of $t_{i_j}$'s with total degree $(k-2)$ 
by Claim 2, and by induction,
$t'$ has the desired form, and so does $t$ by Claim 2.
Otherwise, $a(bc)=\pm (t_1t_2)t_3$.
Consider a next part $s$ so that $t$ has a part $s((t_1t_2)t_3)$ or $((t_1t_2)t_3)s$.
($s$ can be one of $t_1,t_2,t_3$.)
By induction and Claim 1 and 2, $s=\pm t_iz_1$, $\pm t_it_jz_2$ ($i\neq j$)
or $\pm (t_1t_2)t_3z_3$, where $z_1,z_2,z_3$
are products of even power of $t_{i_j}$'s and so $z_1,z_2,z_3\in Z$.
If $z_i\neq 1$, then one can use induction again for the part after taking off $z_i$ from $t$.
Hence we can assume that 
$s=t_i$, $t_it_j$ ($i\neq j$)
or $\pm (t_1t_2)t_3$.
Then $s((t_1t_2)t_3)$ or $((t_1t_2)t_3)s$ has, correspondingly, the factor $t_i^2$, $(t_it_j)^2=-t_i^2t_j^2$ or 
$((t_1t_2)t_3)^2=t_1^2t_2^2t_3^2$ (by Claim 1 and 3),
and so $t=\pm t't_i^2$, $t=\pm t't_i^2t_j^2$ or $t=\pm t't_1^2t_2^2t_3^2$
for some $t'$ which is a product of $t_{i_j}$'s with total degree $k-2$, $k-4$ or $k-6$.
Hence, by induction, $t'$ has the desired form, and so does $t$ by Claim 2.
The case (ii) returns to the case (i) by Claim 1.
For (iii), two of $a,b,c,d$ should be the same,
and so by Claim 1, 2 and the middle Moufang identity,
$(ab)(cd)=\pm t_p^2(t_qt_r)$ for some $p,q,r\in\{1,2,3\}$.
Hence, by the same argument in the first case of (i),
$t$ has the desired form.
\qed
\enddemo

Now, let 
$\Phi [z_1,z_2,z_3]$ be the ordinary polynomial algebra over $\Phi$
in three variables, and let
$$D=(\Phi [z_1,z_2,z_3],z_1,z_2,z_3)$$ be the octonion algebra,
i.e.,
the Cayley-Dickson process over $\Phi [z_1,z_2,z_3]$
three times with structure constants
$z_1$, $z_2$ and $z_3$ starting with trivial involution.
Let $v_1$, $v_2$ and $v_3$ be the basic generators in each step of $D$
so that $v_1^2=z_1$, $v_2^2=z_2$ and $v_3^2=z_3$.
Then $D$ 
has a natural $\Bbb N^3$-grading, 
defining $\deg v_1=(1,0,0)$, $\deg v_2=(0,1,0)$ and $\deg v_3=(0,0,1)$.
It is easily seen that every homogeneous space is a 1-dimensional
free $\Phi $-module.

\example{Remark 3.2}
The algebra $D$ appears as a subalgebra of a free alternative algebra,
discovered by Dorofeev (see  [SSSZ, Theorem 13,  p.296]).
More precisely, let $\Cal F$ be the free alternative algebra over $\Phi $
generated by distinct elements $a$, $b$ and $c$.
Let $u=[a,b]$, $v=(a,b,c)$ and $w=(u,v,a)$.
Then the subalgebra of $\Cal F$ generated by $u$, $v$ and $w$
is isomorphic to $D$
via $u\mapsto v_1$, $v\mapsto v_2$ and $w\mapsto v_3$.
\endexample

\medskip

We now prove our main theorem.

\proclaim{Theorem 3.3}
$\Phi _C[t_1,t_2,t_3]$ is graded isomorphic to $D$.
In particular, 
$Z=\Phi [t_1^2,t_2^2,t_3^2]$, which is the ordinary 
polynomial algebra over $\Phi$ in three variables
$t_1^2$, $t_2^2$, $t_3^2$, and
$\Phi _C[t_1,t_2,t_3]$ is an octonion algebra
$(Z, t_1^2,t_2^2,t_3^2)$.
\endproclaim
\demo{Proof}
By Lemma 1.2,
there exists the epimorphism from 
$\Phi _C[t_1,t_2,t_3]$ onto $D$ defined by
$t_1\mapsto v_1$, $t_2\mapsto v_2$ and $t_3\mapsto v_3$.
So it is enough to show that 
every homogeneous space for the natural $\Bbb N^3$-grading of $\Phi _C[t_1,t_2,t_3]$ 
is generated by one element.
But this follows from Proposition 3.1.
\qed
\enddemo

We note that 
the multiplicative subset 
$$S:=\{z_1^pz_2^qz_3^r\}_{p,q,r\in \Bbb N}$$
of the center $\Phi [z_1,z_2,z_3]$
of $D$
does not contain zero divisors of 
the octonion algebra $D$
(which is 8-dimensional over the center),
and so
the ring of quotients $S^{-1}D$
is also 8-dimensional over the center $\Phi [z_1^{\pm 1},z_2^{\pm 1},z_3^{\pm 1}]$
with generators $v_1$, $v_2$ and $v_3$,
and the multiplication table respect to the generators is the same as the 
multiplication table on $D$.
Hence it
is the octonion algebra 
$$S^{-1}D=(\Phi [z_1^{\pm 1},z_2^{\pm 1},z_3^{\pm 1}], z_1, z_2,z_3).$$
Note that $v_1$, 
$v_2$
and
$v_3$ are invertible in $S^{-1}D$,
and $v_1^{-1}=z_1^{-1}v_1$, 
$v_2^{-1}=z_2^{-1}v_2$
and
$v_3^{-1}=z_3^{-1}v_3$.
Thus 
defining 
$\deg v_1^{-1}=(-1,0,0)$, $\deg v_2^{-1}=(0,-1,0)$ and $\deg v_3^{-1}=(0,0,-1)$,
$S^{-1}D$ has a $\Bbb Z^3$-grading,
and the $\Bbb Z^3$-graded algebra $S^{-1}D$ is called
the {\it octonion $3$-torus}
or the {\it Cayley torus}.
Note that $D$ embeds into $S^{-1}D$.
Clearly, 
every homogeneous space of the Cayley torus is a 1-dimensional
free $\Phi $-module.

\proclaim{Corollary 3.4}
Let $T=\{t_1^pt_2^qt_3^r\}_{p,q,r\in 2\Bbb N}$ be the subset
of $\Phi _C[t_1,t_2,t_3]$.

(1) $T$ is a multiplicative subset of the center of $\Phi _C[t_1,t_2,t_3]$,
which does not contain zero divisors of  $\Phi _C[t_1,t_2,t_3]$,
and $T^{-1}\Phi _C[t_1,t_2,t_3]$
%=(\Phi [t_1^{\pm 2},t_2^{\pm 2},t_3^{\pm 2}], t_1^2,t_2^2,t_3^2)$, which 
can be identified with the Cayley torus
$S^{-1}D$
via $t_1\mapsto v_1$, $t_2\mapsto v_2$ and $t_3\mapsto v_3$.

(2) If $A$ is an alternative algebra over $\Phi $
generated by $a_1$, $a_2$ and $a_3$,
and they satisfy the Cayley relations $a_1a_2=-a_2a_1$,
$a_1a_3=-a_3a_1$, $a_2a_3=-a_3a_2$ and $(a_1a_2)a_3=-a_1(a_2a_3)$,
then $A$ is a homomorphic image of 
$\Phi _C[t_1,t_2,t_3]$
via $t_1\mapsto a_1$, $t_2\mapsto a_2$ and $t_3\mapsto a_3$,
and $a_1^2$, $a_2^2$ and $a_3^2$ are central in $A$.

If, moreover, $a_1$, $a_2$ and $a_3$ are invertible,
then $A$ is a homomorphic image of 
the Cayley torus $T^{-1}\Phi _C[t_1,t_2,t_3]$
via the same map.

(3) The Cayley torus has a presentation
in the category of alternative algebras;
generators
$t_1^{\pm 1}$, $t_2^{\pm 1}$ and $t_3^{\pm 1}$
with relations
$t_it_i^{-1}=1$  for $i=1,2,3$ and 
the Cayley relations (C).
\endproclaim
\demo{Proof}
(1) is now clear by Theorem 3.3.
For (2), let $\vvp:\Phi _C[t_1,t_2,t_3]\longrightarrow A$ be the epimorphism
defined by $t_1\mapsto a_1$, $t_2\mapsto a_2$ and $t_3\mapsto a_3$.
Then the elements of $\vvp(T)$ are central in $A$ by Theorem 3.3.
For the second statement, since $\vvp(T)$ are invertible,
 $\vvp$ extends to 
 $T^{-1}\Phi _C[t_1,t_2,t_3]$
 by the universal property of the ring of quotients. 
 
 For (3), let $Q$ be the alternative algebra having the presentation in the assertion.
 Define $\deg t_1^{\pm 1}=(\pm 1, 0,0)$,
 $\deg t_2^{\pm 1}=(0,\pm 1,0)$ and
 $\deg t_3^{\pm 1}=(0,0\pm 1)$.
 Let $Q^{\alpha}$ be the space generated by the monomials 
 of degree $\alpha\in\Bbb Z^3$.
 Then $Q=\sum_{\alpha\in\Bbb Z^3}\ Q^\alpha$.
 Since the Cayley torus $T^{-1}\Phi _C[t_1,t_2,t_3]$
 has the relations in the assertion, there is a natural homomorphism from $Q$
 onto the Cayley torus
 so that $Q^\alpha$ is mapped onto the 
 homogeneous space of degree  $\alpha$ in the Cayley torus.
 Hence $Q=\oplus_{\alpha\in\Bbb Z^3}\ Q^\alpha$
  (the sum becomes direct).
 On the other hand, by (2), there is a natural graded homomorphism
 from $T^{-1}\Phi _C[t_1,t_2,t_3]$ onto $Q$.
 Hence they are graded isomorphisms.
 \qed
\enddemo

We note that Part (3) of Corollary 3.4 was obtained independently by Bruce
Allison (unpublished).
Because of the presentation of the Cayley torus,
it is reasonable to write
$$T^{-1}\Phi _C[t_1,t_2,t_3]=\Phi _C[t_1^{\pm 1},t_2^{\pm 1},t_3^{\pm 1}].$$
Also, an {\it octonion $n$-torus} ($n\geq 3$) 
$(\Phi [z_1^{\pm 1},\ldots,z_n^{\pm 1}], z_1, z_2,z_3)$
can be written as
$$
\align
\Phi _C[t_1^{\pm 1},\ldots, t_n^{\pm 1}]:
&=\Phi _C[t_1^{\pm 1},t_2^{\pm 1},t_3^{\pm 1}]\otimes_\Phi \Phi[t_4^{\pm 1},\ldots,t_n^{\pm 1}]\\
\endalign
$$
where $\Phi[t_4^{\pm 1},\ldots,t_n^{\pm 1}]$ is the Laurent polynomial algebra over $\Phi$
in $(n-3)$-variables. 
Letting $P:=\Phi[t_4^{\pm 1},\ldots,t_n^{\pm 1}]$,
the octonion torus $\Phi _C[t_1^{\pm 1},\ldots, t_n^{\pm 1}]$ 
can be considered as the Cayley torus over $P$, i.e.,
$P_C[t_1^{\pm 1},t_2^{\pm 1},t_3^{\pm 1}]$.

\medskip

One can start with the alternative algebra $\Phi _C[t_1,\ldots,t_n]$
 over $\Phi$
with generators 
$t_1$, $\ldots$, $t_n$ ($n\geq 3$) and the 
Cayley relations (C),
and the central relations
$$
\text{$[t_i,t_k]=(t_i,t_j,t_k)=0$ 
for $i<j<k$ with $i,j=1,\ldots,n$ and $k=4,\ldots,n$}.
\tag Z
$$
Then by the same argument as above, 
we obtain the following:
\proclaim{Theorem 3.5}
$\Phi _C[t_1,\ldots,t_n]$ is graded isomorphic to $(\Phi [z_1,\ldots,z_n], z_1, z_2,z_3)$.
In particular, the center 
$Z=\Phi[t_1^{2},t_2^{2},t_3^{2},t_4,\ldots,t_n]$, and
$\Phi _C[t_1,\ldots,t_n]$ is an octonion algebra
$(Z, t_1^2,t_2^2,t_3^2)$.

Moreover, generators
$t_1^{\pm 1},\ldots,t_n^{\pm 1}$
with the Cayley relations (C), the central relations (Z),
and the invertible relations $t_it_i^{-1}=1$ for $i=1,\ldots,n$
give a presentation of an octonion torus 
$\Phi _C[t_1^{\pm 1},\ldots, t_n^{\pm 1}]$ in the category of alternative algebras.
\endproclaim

An octonion torus (under the name of the alternative torus)
was found in [BGKN] 
on the classification of extended affine Lie algebras
(under the name of quasi-simple Lie algebras).
The generators and relations of an octonion torus  will be useful for determining generators and relations
for certain extended affine Lie algebras.

\medskip

The following theorem gives a presentation of an octonion algebra over $\Phi $.
\proclaim{Theorem 3.6}
Any octonion algebra over $\Phi $, say $(\Phi ,\mu_1,\mu_2,\mu_3)$
for cancellable scalars $\mu_1,\mu_2,\mu_3$ of $\Phi $,
is isomorphic to $\Phi _C[t_1,t_2,t_3]/(t_1^2-\mu_1, t_2^2-\mu_2, t_3^2-\mu_3)$.
Hence $(\Phi ,\mu_1,\mu_2,\mu_3)$ has a presentation
in the category of alternative algebras;
generators $t_1$, $t_2$ and $t_3$
with Cayley relations (C) and
$t_1^2=\mu_1$, $t_2^2=\mu_2$ and $t_3^2=\mu_3$.

In particular, if $\Phi$ is a field, $A$ is an alternative algebra over $\Phi $
generated by $a_1$, $a_2$ and $a_3$,
and they satisfy $a_1^2=\mu_1$, $a_2^2=\mu_2$, $a_3^2=\mu_3$,
and the Cayley relations $a_1a_2=-a_2a_1$,
$a_1a_3=-a_3a_1$, $a_2a_3=-a_3a_2$ and $(a_1a_2)a_3=-a_1(a_2a_3)$,
then $A$ is isomorphic to $(\Phi ,\mu_1,\mu_2,\mu_3)$.

\endproclaim
\demo{Proof}
Let $B:=\Phi _C[t_1,t_2,t_3]/(t_1^2-\mu_1, t_2^2-\mu_2, t_3^2-\mu_3)$.
Let $v_1$, $v_2$ and $v_3$ be the basic generators of $(\Phi ,\mu_1,\mu_2,\mu_3)$
so that $v_1^2=\mu_1$, $v_2^2=\mu_2$ and $v_3^2=\mu_3$.
Let $\vvp: \Phi _C[t_1,t_2,t_3]
\longrightarrow (\Phi ,\mu_1,\mu_2,\mu_3)$ be the epimorphism
defined by $t_1\mapsto v_1$, $t_2\mapsto v_2$ and $t_3\mapsto v_3$.
Since the ideal $(t_1^2-\mu_1, t_2^2-\mu_2, t_3^2-\mu_3)$ is contained in the kernel of $\vvp$,
$\vvp$ descends to an epimorphism $\ol\vvp:B\longrightarrow (\Phi ,\mu_1,\mu_2,\mu_3)$. 
By Theorem 3.3, 
$\Phi _C[t_1,t_2,t_3]$ is an 8-dimensional free $\Phi [t_1^2,t_2^2,t_3^2]$-module
with basis
$\{1, t_1, t_2, t_3, t_1t_2, t_1t_3, t_2t_3,(t_1t_2)t_3\}$.
So $\{1, \ol t_1, \ol t_2, \ol t_3,\ol{t_1t_2},\ol{t_1t_3},\ol{t_2t_3},\ol{(t_1t_2)t_3}\}$ 
generates $B$ over $\Phi$, where
$\ol{\phantom a}$ is the canonical map from
$\Phi _C[t_1,t_2,t_3]$ onto $B$.
Since $\ol\vvp(\ol t_i)=v_i$ ($i=1,2,3$)
and $\{1, v_1, v_2, v_3, v_1v_2, v_1v_3, v_2v_3,(v_1v_2)v_3\}$
is linearly independent over $\Phi$,
$\ol\vvp$ is injective and
we get $B\cong (\Phi ,\mu_1,\mu_2,\mu_3)$.

For the second statement, one gets $A\cong B$
since $A$ is a homomorphic image
of the simple algebra $B$ and $A\neq0$.
\qed
\enddemo

\head
\S4 Cayley-Dickson rings
\endhead

Let $R$ be an alternative ring with a nonzero center $Z$
which does not contain zero divisors of $R$
(e.g. $R$ is prime).
Then $Z^*=Z\setminus\{0\}$ is a multiplicative subset of $Z$, and one can construct
the ring of quotients $(Z^*)^{-1}R$,
which is called the {\it central closure} of $R$,
denoted $\ol R$.
We note that $R$ embeds into $\ol R$, $\ol Z=(Z^*)^{-1}Z$ is a field
of fractions of $Z$,
$\ol R$ is a central $\ol Z$-algebra,
and $\ol R\cong \ol Z\otimes_ZR$.
Moreover, $R$ is called a {\it Cayley-Dickson ring} if
the central closure $\ol R$ is an octonion algebra over $\ol Z$
(see [SSSZ, p.193]).
For example, if $\Phi$ is a domain, our ring of Cayley polynomials $\Phi_C[t_1,t_2,t_3]$ or
an octonion torus $\Phi _C[t_1^{\pm 1},\ldots, t_n^{\pm 1}]$
is a Cayley-Dickson ring
so that the central closure is
$(\ol Z,t_1^2,t_2^2,t_3^2)$,
where $\ol Z=\ol\Phi(t_1^2,t_2^2,t_3^2)$ or $\ol\Phi(t_1^2,t_2^2,t_3^2, t_4,\ldots,t_n)$
(rational function fields over the field of fractions $\ol\Phi$ in 3 or $n$ variables),
respectively.
However, a Cayley-Dickson ring is not necessarily an octonion algebra over the center
(see Example 4.3).

\proclaim{Lemma 4.1}
Let $A$ be an alternative algebra over $\Phi$ with center $Z$
which does not contain zero divisors of $A$.
Assume that $A$ is generated by $a_1$, $a_2$ and $a_3$,
and they satisfy the Cayley relations $a_1a_2=-a_2a_1$,
$a_1a_3=-a_3a_1$, $a_2a_3=-a_3a_2$ and $(a_1a_2)a_3=-a_1(a_2a_3)$,
then $A$ is an octonion algebra over $Z$, isomorphic to $(Z ,a_1^2,a_2^2,a_3^2)$.

\endproclaim
\demo{Proof}
Since $A$ embeds into $\ol A$,
$a_1$, $a_2$ and $a_3$ also satisfy the Cayley relations in $\ol A$.
Also, $a_1^2,a_2^2,a_3^2\in Z$ by Corollary 3.4 (2).
Hence, by Theorem 3.6,
$\ol A$ is an octonion algebra over the field $\ol Z$,
i.e., $\ol A=(\ol Z,a_1^2,a_2^2,a_3^2)$.
In particular, $A$ is an 8-dimensional free $Z$-module
with basis
$\{1, a_1, a_2, a_3, a_1a_2, a_1a_3, a_2a_3,(a_1a_2)a_3\}$,
which is also a basis of $\ol A$.
Moreover, $A$ and $\ol A$ has the same multiplication table
relative to the basis, and so $A=(Z,a_1^2,a_2^2,a_3^2)$.
\qed
\enddemo

Let us state a celebrating theorem in alternative theory [SSSZ, Theorem 9, p.194]:
\proclaim{Slater's Theorem}
Any prime nondegenerate alternative algebra that is not associative
is a Cayley-Dickson ring.
\endproclaim

(It is also true that every Cayley-Dickson ring is a prime nondegenerate ring
[SSSZ, Proposition 3, p.193].) 
Using the theorem, we have the following:

\proclaim{Proposition 4.2}
Let $R$ be a prime nondegenerate alternative algebra that is not associative
over $\Phi$, and $Z=Z(R)$ its center.
Then there exist a subalgebra $A$ of $R$
so that $Z(A)\subset Z$
and $A$ is an octonion algebra over the center $Z(A)$.
Moreover, the subalgebra $B:=ZA$ is an octonion algebra over the center $Z$ and
the central closures of $B$ and $R$ coincide, i.e.,
$\ol B=\ol {R}$.

Also, $\ol R$ is a base field extension of $\ol A$, namely,
$\ol R\cong\ol Z\otimes_K\ol A$, where $K=\ol {Z(A)}$.
\endproclaim
\demo{Proof}
By Slater's Theorem, 
$\ol R$ is an octonion algebra over the field $\ol Z$.
Let $v_1$, $v_2$ and $v_3$ be the basic generators and so
they satisfy the Cayley relations.
Note that $v_1=z_1^{-1}r_1$, $v_2=z_2^{-1}r_2$ and $v_3=z_3^{-1}r_3$
for some $z_1,z_2,z_3\in Z^*$ and $r_1,r_2,r_3\in R$.
So $a_1:=z_1z_2z_3v_1$, $a_2:=z_1z_2z_3v_2$ and $a_3:=z_1z_2z_3v_3$
also satisfy the Cayley relations and they are in $R$.
Let $A$ be the subalgebra of $R$ generated by $a_1$, $a_2$ and $a_3$.
Note that 
if $z\in Z(A)$, then $z$ is, in particular, central for a generating set $\{a_1,a_2,a_3\}$ of $A$,
and so is for a generating  set $\{v_1,v_2,v_3\}$ of $\ol R$.
Hence, by the theorem of Bruck and Kleinfeld [SSSZ, Lemma 16, p.289], 
$z$ is central for $\ol R$, and so $z\in Z$.
Thus, $Z(A)$ does not contain zero divisors of  $A$, and
hence
by Lemma 4.1, $A=(Z(A),a_1^2,a_2^2,a_3^2)$.

Now, we have $Z\subset Z(B)\subset Z$,
and so $Z=Z(B)$.
Thus by Lemma 4.1 again,
$B=(Z,a_1^2,a_2^2,a_3^2)$.
Finally, for any $r\in R$,
there exists some $z\in Z^*$ such that $zr\in B$.
In fact, 
$r=f(v_1,v_2,v_3)=f(z_1^{-1}z_2^{-1}z_3^{-1}a_1,z_1^{-1}z_2^{-1}z_3^{-1}a_2,z_1^{-1}z_2^{-1}z_3^{-1}a_3)$ 
for some polynomial $f$ over $\ol Z$.
So there exists $z\in Z^*$ such that 
$zf(z_1^{-1}z_2^{-1}z_3^{-1}a_1,z_1^{-1}z_2^{-1}z_3^{-1}a_2,z_1^{-1}z_2^{-1}z_3^{-1}a_3)=g(a_1,a_2,a_3)$
for some polynomial $g$ over $Z$. Hence $zr=g(a_1,a_2,a_3)\in B$.
Thus $r=z^{-1}b$ for some $b\in B$, and so $\ol R\subset\ol{B}$.
Since the other inclusion is clear,
we obtain $\ol R=\ol{B}$.

For the last statement, let $\vvp$ be a $K$-linear map from 
$\ol Z\otimes_K\ol A$ to $\ol R$ defined by $\vvp(u_i\otimes w_j)=u_iw_j$ 
for a basis $\{u_i\}$ of $\ol Z$ over $K$
and a basis $\{w_i\}$ of $\ol A$ over $K$.
%(Note that $\ol Z\otimes_K\ol A=\ol Z\otimes_K(K\otimes_{Z(A)}A)=\ol Z\otimes_{Z(A)}A$.)
Then $\vvp$ is a homomorphism, and $\ol Z$-linear. 
Since $\ol Z\otimes_K\ol A$ and $\ol R$ are both 8-dimensional over $\ol Z$,
it is enough to show that $\vvp$ is onto.
For $v^{-1}r\in\ol R$ ($v\in Z^*$ and $r\in R$),
there exists $z\in Z^*$ and $b\in B$ such that $zr=b$  by the above,
and so $zr=\sum_k\ z_ka_k$  for some $z_k\in Z$ and $a_k\in A$.
Hence $v^{-1}r=v^{-1}z^{-1}\sum_k\ z_ka_k=\vvp(\sum_k\ v^{-1}z^{-1}z_k\otimes a_k)$, and so $\vvp$ is onto.
\qed
\enddemo

Let us finally give an example of a prime nondegenerate algebra which is not an octonion algebra
over the center.
\example{Example 4.3}
For simplicity, let $F$ be a field of characteristic $\neq 2$,
and let $F[z]$ be the ordinary polynomial algebra over $F$.
Let $F[z]_C[t_1,t_2,t_3]$ be the ring of Cayley polynomials over $F[z]$.
Let $R$ be the $F$-subalgebra of $F[z]_C[t_1,t_2,t_3]$ generated by
$t_1$, $t_2$, $t_3$ and $zt_1$.
Then the center $Z=Z(R)=F[t_1^2,t_2^2,t_3^2,z^2t_1^2,zt_1^2]$,
and $R$ is a 12-dimensional free $Z$-module with basis
$$
\{1, t_1, t_2, t_3, t_1t_2, t_1t_3, t_2t_3, (t_1t_2)t_3,
zt_1, zt_1t_2, zt_1t_3,   zt_1(t_2t_3)
\}.
$$
Hence $R$ is not an octonion algebra over $Z$.
But the central closure $\ol R$ is an octonion algebra over
$\ol Z=F(z, t_1^2,t_2^2,t_3^2)$,
i.e.,
$\ol R=(\ol Z, t_1^2,t_2^2,t_3^2)$
(by Theorem 3.6),
and so $R$ is a Cayley-Dickson ring.
One can take a subalgebra $A$ of $R$ in Proposition 4.2 as
$A=F_C[t_1,t_2,t_3]$.
\endexample

%\enskip

\end